\documentclass[12pt,a4paper]{article}

\usepackage{amssymb,amsmath,amsfonts,amsthm,enumerate}

\usepackage{graphicx} %for figures
\usepackage{float} %for exact placement of figure (H)

\usepackage{tikz}

\begin{document}

\centerline{\bf{Existence results for anisotropic and isotropic $p(x)$-Laplace equations}}

\bigskip

\centerline{Alkis S. Tersenov} 

\bigskip

\begin{small}
\centerline {$^a$University of Crete, Department of Mathematics and Applied Mathematics,} 
\centerline{71003 Heraklion -- Crete, Greece;}
\centerline {e-mail: tersenov@uoc.gr}

\bigskip

\noindent
\textbf{Abstract.} The Dirichlet problem is considered both for degenerate and singular inhomogeneous quasilinear parabolic equations. We prove the existence of a solution $u$ such that $u_t$ belongs to 
$L_{\infty}$. The $L_{\infty}$ estimate of $u_t$ is obtained by introducing a new time variable. 

\medskip


\medskip

\end{small}

\medskip

\bigskip

\centerline{{\bf \S 1. Introduction and Main Results}}

\medskip

Let $\Omega$ be a bounded domain in ${\bf R}^n$ satisfying the exterior sphere condition and
$Q_T=(0,T)\times \Omega$ with an arbitrary $T\in(0,\infty)$. By $x=(x_1,\dots,x_n)$ we denote points in 
$\Omega$ and by $t$ the time variable that varies in the interval $[0,T]$.

1. \underline{Anisotropic case}.

Consider the following anisotropic equation  
$$
u_t=\sum_{i=1}^{n}(|u_{x_i}|^{p_i(x)} u_{x_i})_{x_i}+F \ \ 
\hbox{in}\ \ Q_T,\ \ \min_{\overline{\Omega}}p_i(x)> -1, \leqno (1.1)
$$
coupled with the Dirichlet boundary condition
$$
u=\psi(x)\ \ \hbox{on}\ \ \partial\Omega\times (0,T) \leqno (1.2)
$$
and the initial condition
$$
u(0,x)=u_0(x)\ \ \hbox{in}\ \ \Omega\ \ \hbox{where}\ \ u_0(x)=\psi(x) \ \ \hbox{on}\ \ \partial\Omega. \leqno (1.3)
$$
Assume that 
$$
F=\sum_{i=1}^{n}f_i(x)u_{x_i} +f_0(x,u), \leqno (1.4)
$$
where   
$$
f_0(x,u)-f_0(x,v)\leq 0\ \ \hbox{if}\ \ u> v, \leqno (1.5)
$$
and
$$
uf_0(x,u)\leq a_1u^2+a_2. \leqno (1.6)
$$
Here $a_1$, $a_2$ are some positive constants. Condition (1.6) guaranties (see Section 2) the a priori estimate 
$$
|u|\leq M
$$
with $M$ depending only on $\max |u_0|$ and $a_1,\, a_2$. 

Concerning the initial function $u_0$ we assume that
$$
\max \Big|\sum_{i=1}^{n}(|u_{0x_i}|^{p_i(x)} u_{0x_i})_{x_i}
+F(x,u,\nabla u_0)\Big|<+\infty \leqno (1.7)
$$
where maximum is taken over the set  $\overline{\Omega}\times [-M-M]$. The functions 
$f_i$ and $p_i$ obey the following condition
$$
p_i(x)\in C^{1+\gamma}(\overline{\Omega}), f_i(x)\in C^{0}(\overline{\Omega})\, i=1,...,n,\ \ \, 
f_0(x,u)\in C^{0}(\overline{\Omega}\times[-M,M]) \leqno (1.8)
$$
for some $\gamma\in(0,1)$.

\noindent
{\bf Definition 1.} {\it We say that a function $u$ is a solution of problem 
$(1.1)$ - $(1.3)$ if
$$
u \in L_{\infty}(Q_T),\ \ u_{x_i}\in L_{\infty}(0,T;L_{p_i+2}(\Omega)),\ \ u_t\in L_{\infty}(Q_T)  
$$
and the following integral identity
$$
\int_{Q_T}\Big(u_t\,\phi+\sum_{i=1}^n|u_{x_i}|^{p_i}u_{x_i}\,\phi_{x_i}\Big)dxdt=\int_{Q_T} F\phi dxdt
$$
holds for an arbitrary smooth function $\phi$ which vanishes on $(0,T)\times \partial\Omega$. The initial condition is satisfied in the classical sense and the boundary condition in the sense of trace.
}

\medskip

\noindent
{\bf Theorem 1.} {\it Suppose that $(1.4)$ - $(1.8)$ are fulfilled, then for an arbitrary $T>0$, 
there exists a unique solution of problem $(1.1)$ - $(1.3)$. 
}
 
\medskip

2. \underline{Isotropic case}.

Consider the following non homogeneous quasilinear parabolic equation  
$$
u_t=div \big(|\nabla u|^{p(x)} \nabla u\big)+F \ \ \hbox{in}\ \ Q_T,\ \ 
\min_{\overline{\Omega}}p(x)>-1, \leqno (1.9)
$$
coupled with conditions (1.2), (1.3). Here $\nabla u=(u_{x_1},...,u_{x_n})$. Concerning the function $u_0$ we assume that
$$
\max \Big|div (|\nabla u_0|^{p(x)} \nabla u_0)
+F(x,u,\nabla u_0)\Big|<+\infty \leqno (1.10)
$$
where maximum is taken over the set  $\overline{\Omega}\times [-M-M]$. 

\smallskip

\noindent
{\bf Definition 2.} {\it We say that a function $u$ is a solution of problem 
$(1.9)$, $1.2$, $(1.3)$ if
$$
u \in L_{\infty}(Q_T),\ \ |\nabla u|\in L_{\infty}(0,T;L_{p+2}(\Omega)),\ \ u_t\in L_{\infty}(Q_T)  
$$
and the following integral identity
$$
\int_{Q_T}\Big(u_t\,\phi+|\nabla u|^{p}\nabla u\cdot \nabla\phi\Big)dxdt=\int_{Q_T} F\phi dxdt
$$
holds for an arbitrary smooth function $\phi$ which vanishes on $(0,T)\times \partial\Omega$. The initial condition is satisfied in the classical sense and the boundary condition in the sense of trace.
}

\smallskip

\noindent
{\bf Theorem 2.} {\it Suppose that $(1.4)$ - $(1.6)$, $(1.8)$ and $(1.10)$ are fulfilled, then for an arbitrary $T>0$, there exists a unique weak solution of problem $(1.9)$, $(1.2)$, $(1.3)$. 
}

\medskip

In the last decade, equation (1.1) with non-constant exponents $p_i$ has attracted the attention of many researchers, see, for example, [2], [3], [9], [20], [27]. Concerning the physical motivation see [2] and the references therein. The existence of a weak solution of the Dirichlet problem 
belonging to $L_{\infty}(0,T;L_2(\Omega))$ with first derivative $u_{x_i}$ 
from $L_{p_i+2}(Q_T)$ and $u_t$ belonging to the corresponding negative space was proved in [2] 
where a more general than (1.1) equation was considered. Equation (1.9) has been intensively studied in recent years by many authors. See [4], [5] - [7], [10], [11], [14], [16], [17], [20], [25], [28], [29], where the time derivative belongs to $L_2(Q_T)$ or to some negative space. In [16] the following equation was considered 
$$
u_t- div (|\nabla u|^{p(x,t)}\nabla u) = f(t,x)\ \ \hbox{in}\ \  Q_T
$$
with smooth $\partial\Omega$ and 
$$
p(t,x)>\frac{-4}{n+2}.
$$
It was shown that if $p, f$ are H\"{o}lder continuous in $Q_T$ and $|\nabla u_0|^{p(0,x)}\in L_1(\Omega),$ 
then $u\in C_{t,x}^{1/2,1}(Q_T\cap \{t\ge\epsilon\})$ for every $\epsilon > 0$. 
Moreover, if $f$ is Lipschitz-continuous in $t$,
$p =p(x),$ $\nabla p\in L_{\infty}(\Omega)$ and either $p(x)>0$, or 
$\frac{-4}{n+2}<p(x)< 0,$ then the solutions are Lipschitz continuous in time in $Q_T\cap \{t\ge\epsilon\}$. Obviously, for $n=1$ equations (1.1) and (1.9) coincide. The one dimensional case with zero boundary conditions and $p=p(x)$ was investigated in [20]. In the present paper, we extend the results obtained in [20] to the multidimensional case with non homogeneous boundary conditions in two directions: anisotropic (1.1) and isotropic (1.9). As far as we know, for $n>1$ the boundedness of $u_t$ for both anisotropic and isotropic cases is a new result. Note that for a multidimensional anisotropic homogeneous equation with zero boundary conditions and constant exponents $p_i$, the $L_{\infty}$ estimate of $u_t$ was obtained in [22]. 

\medskip

In order to prove Theorems 1 and 2  we regularize our equation and obtain a priori estimates (independent from regularization) which allow us to pass to the limit in the regularized problem and to obtain the needed solution. The key estimate is the $L_{\infty}$ estimate of $u_t$ which is obtained by introducing a new time variable. This method was first proposed in [20] by analogy with Kruzhkov's method of introducing a new spatial  variable (see, for example, [21]). 

Note that the a priori estimates for the regularized problem will be obtained under more general right hand side than (1.4). See the Remark at the very end of the last section on this paper. 
 
\bigskip

\smallskip

\centerline {{\bf \S 2. Regularized Problem}}

\medskip

Without loss of generality we assume that $f_0(x,u)$ is continuously differentiable functions of its arguments. 
Otherwise instead of $f_0$ we take 
$$
f_{0\varepsilon}(x,u)=\frac{1}{2\varepsilon}
\int_{\Omega}\Big(\int_{u-\varepsilon}^{u+\varepsilon}f_0(y,\zeta)d\zeta\Big)\eta_{\varepsilon}(|x-y|)dy,
$$ 
here $\eta_{\varepsilon}$ is a standard Steklov average (mollifier) (see [13]). Obviously 
$\frac{\partial}{\partial u}f_0(\varepsilon,x,u)\leq 0$ and hence condition (1.5) (as well as (1.6) ) 
for $f_0(\varepsilon,x,u)$ is satisfied. Obviously, $f_{0\varepsilon}\to f_0$ a.e. as $\varepsilon\to 0$,
moreover, since $f_0(x,u)$ is continuous, uniform convergence takes place. 

Similarly, without loss of generality we assume that $f_i(x),\,p_{ix_i}(x)$ ($i=1,..., n$) and $p_{x_i}(x)$ are continuously differentiable functions. 

\smallskip

1. \underline{Anisotropic case}.

Regularize equation (1.1):
$$
u_{\varepsilon t}=\sum_{i=1}^n((u_{\varepsilon {x_i}}^{2}
+\varepsilon)^{p_i/2}u_{\varepsilon {x_i}})_{x_i}+F(x,u,\nabla u) \ \ \hbox{in}\ \ Q_T \leqno (2.1)
$$
where $\varepsilon\in (0, \varepsilon_0]$ for some strictly positive constant $\varepsilon_0$. Rewrite this equation in the equivalent form 
$$
u_{\varepsilon t}=\sum_{i=1}^n a_i(\varepsilon,x,u_{\varepsilon x_i})u_{\varepsilon x_ix_i}+
\sum_{i=1}^{n}b_i(\varepsilon,x,u_{\varepsilon x_i})+F(x,u_{\varepsilon},\nabla u_{\varepsilon}), \ \ \hbox{in}\ \ Q_T \leqno (2.2)
$$
where
$$
a_i(\varepsilon,x,z_i)=(z_i^{2}+\varepsilon)^{\frac{p_i-2}{2}}\big((1+p_i)z_i^{2}+\varepsilon\big)
$$
and 
$$
b_i(\varepsilon,x,z_i)=p_{x_i}z_i(z_i^{2}+\varepsilon)^{\frac{p_i}{2}}\ln (z_i^{2}+\varepsilon)^{\frac{1}{2}} 
$$

Concerning the existence of a classical solution $u_\varepsilon$ of problem (2.2), (1.2), (1.3) with $F$ as in (1.4) see [13, Chapter VI]. 

\smallskip 

2. \underline{Isotropic case}.

Regularize equation (1.9):
$$
u_{\varepsilon t}=div \Big((|\nabla u_{\varepsilon}|^2+\varepsilon)^{p(x)/2}\nabla u_{\varepsilon}\Big)+
F(x,u_{\varepsilon},\nabla u_{\varepsilon}) \ \ \hbox{in}\ \ Q_T \leqno (2.3)
$$
where $\varepsilon$ is strictly positive constant. Rewrite this equation in the equivalent form 
$$
u_{\varepsilon t}=\sum_{i=1}^n a_{ij}(\varepsilon,x,\nabla u_{\varepsilon})u_{\varepsilon x_ix_i}+
\sum_{i=1}^{n}b_{i}(\varepsilon, x,\nabla u_{\varepsilon})+F(x,u_{\varepsilon},\nabla u_{\varepsilon}), \ \ \hbox{in}\ \ Q_T \leqno (2.4)
$$
where
$$
a_{ii}(\varepsilon,x,z)=\big(|z|^2+\varepsilon\big)^{\frac{p-2}{2}}
\big(|z|^2+pz_i^2+\varepsilon\big),\ \ z=(z_1,...,z_n)
$$
$$
a_{ij}(\varepsilon,x,z)=p\big(|z|^2+\varepsilon\big)^{\frac{p-2}{2}}
z_iz_j\ \ \hbox{for}\ \ i\not=j
$$
and 
$$
b_{i}(\varepsilon,x,z)=p_{x_i}z_i(|z|^2+\varepsilon)^{\frac{p}{2}}\ln (|z|^2+\varepsilon)^{\frac{1}{2}}. 
$$

Similarly to the anisotropic case the existence of a classical solution 
$u_\varepsilon$ of problem (2.2), (1.2), (1.3) with $F$ as in (1.4) is known (see  [13] (Chapter VI). 

\smallskip

Our goal is to obtain uniform with respect to $\varepsilon$ estimates on $u_{\varepsilon}$ which would allow us
to pass to the limit as $\varepsilon\to 0$. In the next section, we will obtain the $L_{\infty}$ estimate of $u_{\varepsilon t}$ for classical solutions. 

\smallskip

\bigskip

\centerline {{\bf \S 3. A priori estimates of $u_{\varepsilon t}$}}

\medskip

Let us start with the formulation of the standard estimate of $|u_{\varepsilon}|$ (see [13, Chapter 1, \S2 ]. 
Assume that
$$
uF(x,u,0)\leq a_1u^2+a_2  \leqno (3.1)
$$  
for some constant $a_1$, $a_2$.

\noindent
{\bf Lemma 3.1}. {\it If condition $(3.1)$ is fulfilled, then the following estimate takes place 
$$
|u_{\varepsilon}(x,t)|\le M
$$
with $M$ depending only on $\max |u_0|$ and $a_1,\, a_2$.
}

We will obtain an $L_{\infty}$ estimate of  $u_{\varepsilon t}$ for a bounded $u_{\varepsilon}$ with the only restrictions on $F(x,u,\nabla u)$ the following monotonicity assumption and (3.1) 
$$
F(x,u,z)-F(x,v,z)\leq 0\ \ \hbox{for} \ \ u\ge v. \leqno (3.2)
$$

1. \underline{Anisotropic case}.

Consider equation (2.1) coupled with conditions (1.2), (1.3). Denote by $K$ the following quantity:
$$
K=\max\Big|\sum_{i=1}^n((u_{0x_i}^2
+\varepsilon)^{p_i/2}u_{0x_i})_{x_i}\Big| +\max |F(x,u,\nabla u_0)| \leqno (3.3)
$$ 
where maximum is taken over the set $(\varepsilon,x,u)\in 
[0,\varepsilon_0]\times\overline{\Omega}\times  [-M,M].$ It is easy to see that $K<+\infty$ due to (1.7).

\smallskip

Below for simplicity in the proofs we will omit the sub index $\varepsilon$.

\smallskip

\par\noindent
{\bf Lemma 3.2}. {\it For every $\varepsilon\in (0,\varepsilon_0]$ 
the following estimate
$$
|u_{\varepsilon}(t,x)-u_0(x)|\leq K\,t,
$$
takes place.
}
\par\noindent 
{\bf Proof.} Introduce the function
$$
h(t)=(K+\delta)t\ \ \hbox{in}\ \ [0,T],
$$
where $\delta>0$. Let us prove the following inequality
$$
u(t,x)-u_0(x)\leq h(t). \leqno (3.4)
$$
Consider the linear operator 
$$
L\equiv \frac{\partial}{\partial t}-\sum_{i=1}^na_i(\varepsilon,x,u_{0x_i})\frac{\partial^2}{\partial x_i^2}.
$$
Define the function $\phi^+\equiv u-[u_0(x)+h(t)]$, obviously
$$
L\phi^+=u_t-
\sum_{i=1}^n a_i(\varepsilon,x,u_{0x_i})u_{x_ix_i}
+\sum_{i=1}^n a_i(\varepsilon,x,u_{0x_i})u_{0x_ix_i}-K-\delta.
$$

Denote by $\Gamma_T$ the parabolic boundary of $Q_T$, i.e. 
$$
\Gamma_T=\partial Q_T\setminus \{(T,x):x\in {\overline \Omega}\setminus \partial \Omega\}. \leqno (3.5)
$$
Suppose that at the point $N\in {\overline Q_T}\setminus\Gamma_T$ the 
function $\phi^+$ attains its maximum, then at this point 
$$
\nabla \phi^+=0\ \ \ \Leftrightarrow \ \ \nabla u=\nabla u_0 \ \ \Rightarrow\ \ a_i(\varepsilon,x,u_{0x_i})=
a_i(\varepsilon,x,u_{x_i})\ \ \Rightarrow
$$
$$
L\phi^+\Big|_N=
\sum_{i=1}^n b_i(\varepsilon,x,u_{\varepsilon x_i})+
F(x,u,\nabla u)+
\sum_{i=1}^n a_i(\varepsilon,x,u_{0x_i})u_{0x_ix_i}-K-
\delta\Big|_{N}=
$$
$$
\sum_{i=1}^{n}b_i(\varepsilon,x,u_{0 x_i})+\sum_{i=1}^n a_i(\varepsilon,x,u_{0x_i})u_{0x_ix_i}+
F(x,u,\nabla u_0)-K-\delta\Big|_{N}=
$$
$$
\sum_{i=1}^n((u_{0x_i}^{2}+\varepsilon)^{p_i/2}u_{0x_i})_{x_i}
+F(x,u,\nabla u_0)-K-\delta\Big|_{N}<0.
$$
which is impossible (the last inequality is due to (3.3) ). 

Consider $\phi^+$ on $\Gamma_T:$ for $x\in\partial\Omega$, $t\in [0,T]$ we have $\phi^+=-h(t)\leq 0$ and  
for $t=0$, $x\in \Omega$ we have $\phi^+=-h(0)=0$.

Thus $\phi^+\leq 0$  in $Q_T$ and (3.4) is proved.

\smallskip 

Let us show now that
$$
u(t,x)-u_0(x)\ge -h(t).\leqno (3.6)
$$
Introduce the function $\phi^-\equiv u-[u_0(x)-h(t)]$. We have
$$
L\phi^-=u_t-
\sum_{i=1}^n a_i(\varepsilon,x,u_{0x_i})u_{x_ix_i}+\sum_{i=1}^n a_i(\varepsilon,x,u_{0x_i})u_{0x_ix_i}+K+\delta 
$$
Suppose that at the point $N_1\in {\overline Q_T}\setminus\Gamma_T$ the 
function $\phi^-$ attains its minimum, then at this point
$$
\nabla \phi^-=0\ \ \Leftrightarrow\ \ \nabla u=\nabla u_0\ \ \Rightarrow\ \ a_i(\varepsilon,x,u_{0x_i})=
a_i(\varepsilon,x,u_{x_i})
$$
and, similarly to the previous case,
$$
L\phi^-\Big|_{N_1}=\sum_{i=1}^n((u_{0x_i}^{2}+\varepsilon)^{p_i/2}u_{0x_i})_{x_i}
+F(x,u,\nabla u_0)+K+\delta\Big|_{N_1}>0,
$$
which is impossible. Consider $\phi^-$ on $\Gamma_T:$
\par\noindent
for $x\in \partial\Omega$, $t\in [0,T]$ we have $\phi^-=h(t)\ge 0$ and 
for $t=0$, $x\in \Omega$ we have $\phi^-=h(0)=0$.

Thus $\phi^-\ge 0$  and (3.6) is proved. From (3.4) and (3.6) we obtain
$$
|u(t,x)-u_0(x)|\leq h(t),
$$
passing to the limit as $\delta\to 0$ we finish the prove of Lemma 3.2. $\square$

\smallskip

\noindent
{\bf Lemma 3.3}. {\it For every $\varepsilon\in (0,\varepsilon_0]$ 
the inequality 
$$
|u_{\varepsilon t}|\leq K
$$
holds.}
\par\noindent 
{\bf Proof.} 
Consider equation (2.2) in two different points $(t,x)$ and $(\tau, x)$:
$$
u_t-\sum_{i=1}^n a_i(\varepsilon,x,u_{x_i})u_{x_ix_i}=\sum_{i=1}^nb_i(\varepsilon,x,u_{x_i})+F(x,u,\nabla u),\ \ u=u(t,x) \leqno(3.7)
$$
$$
u_{\tau}-\sum_{i=1}^n a_i(\varepsilon,x,u_{x_i})u_{x_ix_i}=\sum_{i=1}^nb_i(\varepsilon,x,u_{x_i})+F(x,u,\nabla u),\ \ u=u(\tau,x). \leqno(3.8)
$$
Subtracting (3.8) from (3.7) for $v(t,\tau,x)\equiv u(t,x)-u(\tau,x)$, since 
$$
v_t=u_t(t,x),\ \ v_{\tau}=-u_{\tau}(\tau,x),\ \ v_{x_ix_i}= u_{x_ix_i}(t,x)-u_{x_ix_i}(\tau,x),
$$
we obtain 
$$
v_t+v_{\tau}-\sum_{i=1}^n a_i(\varepsilon,x,u_{x_i}(t,x))v_{x_ix_i}=
$$
$$
\sum_{i=1}^n\big[a_i(\varepsilon,x,u_{x_i}(t,x))-a_i(\varepsilon,x,u_{x_i}(\tau,x)\big]u_{x_ix_i}(\tau,x)+
$$
$$
\sum_{i=1}^n\big[b_i(\varepsilon,x,u_{x_i}(t,x))-b_i(\varepsilon,x,u_{x_i}(\tau,x)\big]+
$$
$$
F(x,u(t,x),\nabla u(t,x))-F(x,u(\tau,x),\nabla u(\tau,x)).
$$
Consider the function
$$
\hbox{w}\equiv v-K(t-\tau)=u(t,x)-u(\tau,x)-K(t-\tau)
$$
in the prism (see Figure \ref{fig1})
$$
P=\{(t,\tau,x):t\in(0,T), \tau\in(0,T),x\in \Omega, t>\tau\}. \leqno(3.9)
$$
%\begin{figure}[H]
%    \centering
%    \includegraphics[width=0.5\textwidth]{prism.eps}
%    \caption{prism $P$  \label{fig1}}
%\end{figure}

\tikzset{every picture/.style={line width=0.75pt}} %set default line width to 0.75pt      
\begin{figure}  
\begin{center}

\begin{tikzpicture}[x=0.75pt,y=0.75pt,yscale=-1,xscale=1]
%uncomment if require: \path (0,261); %set diagram left start at 0, and has height of 261

%Straight Lines [id:da6552435180425036] 
\draw [color={rgb, 255:red, 0; green, 0; blue, 0 }  ,draw opacity=1 ]   (335.33,162.33) -- (472.87,162.33) ;
\draw [shift={(474.87,162.33)}, rotate = 180] [color={rgb, 255:red, 0; green, 0; blue, 0 }  ,draw opacity=1 ][line width=0.75]    (10.93,-3.29) .. controls (6.95,-1.4) and (3.31,-0.3) .. (0,0) .. controls (3.31,0.3) and (6.95,1.4) .. (10.93,3.29)   ;
%Straight Lines [id:da8615690315359394] 
\draw [color={rgb, 255:red, 0; green, 0; blue, 0 }  ,draw opacity=1 ]   (335.33,162.33) -- (335.33,24.47) ;
\draw [shift={(335.33,22.47)}, rotate = 90] [color={rgb, 255:red, 0; green, 0; blue, 0 }  ,draw opacity=1 ][line width=0.75]    (10.93,-3.29) .. controls (6.95,-1.4) and (3.31,-0.3) .. (0,0) .. controls (3.31,0.3) and (6.95,1.4) .. (10.93,3.29)   ;
%Straight Lines [id:da32667350483344304] 
\draw [color={rgb, 255:red, 0; green, 0; blue, 0 }  ,draw opacity=1 ]   (335.33,162.33) -- (227.2,234.03) ;
\draw [shift={(225.53,235.13)}, rotate = 326.45] [color={rgb, 255:red, 0; green, 0; blue, 0 }  ,draw opacity=1 ][line width=0.75]    (10.93,-3.29) .. controls (6.95,-1.4) and (3.31,-0.3) .. (0,0) .. controls (3.31,0.3) and (6.95,1.4) .. (10.93,3.29)   ;
%Straight Lines [id:da14371413515745268] 
\draw [color={rgb, 255:red, 208; green, 2; blue, 27 }  ,draw opacity=1 ][line width=0.75]  [dash pattern={on 4.5pt off 4.5pt}]  (335.33,162) -- (394.2,102.47) ;
%Straight Lines [id:da6267974663888056] 
\draw [color={rgb, 255:red, 208; green, 2; blue, 27 }  ,draw opacity=1 ][line width=0.75]    (257.8,214.07) -- (316.67,154.53) ;
%Straight Lines [id:da6179939624226787] 
\draw [color={rgb, 255:red, 208; green, 2; blue, 27 }  ,draw opacity=1 ][line width=0.75]    (316.67,154.53) -- (394.2,102.47) ;
%Straight Lines [id:da4435489001501298] 
\draw [color={rgb, 255:red, 208; green, 2; blue, 27 }  ,draw opacity=1 ][line width=0.75]    (316.67,154.53) -- (316.2,214.07) ;
%Straight Lines [id:da8504532034640451] 
\draw [color={rgb, 255:red, 208; green, 2; blue, 27 }  ,draw opacity=1 ][line width=0.75]    (394.2,102.7) -- (393.73,162.23) ;
%Straight Lines [id:da9753988576505686] 
\draw [color={rgb, 255:red, 208; green, 2; blue, 27 }  ,draw opacity=1 ][line width=0.75]    (316.67,214.07) -- (394.2,162) ;
%Straight Lines [id:da5858998281266996] 
\draw [color={rgb, 255:red, 208; green, 2; blue, 27 }  ,draw opacity=1 ] [dash pattern={on 4.5pt off 4.5pt}]  (335.33,162) -- (394.2,162) ;
%Straight Lines [id:da3725005995771671] 
\draw [color={rgb, 255:red, 208; green, 2; blue, 27 }  ,draw opacity=1 ]   (257.33,214.07) -- (316.2,214.07) ;
%Straight Lines [id:da12931642265899468] 
\draw [color={rgb, 255:red, 208; green, 2; blue, 27 }  ,draw opacity=1 ][line width=0.75]  [dash pattern={on 4.5pt off 4.5pt}]  (257.8,214.07) -- (335.33,162) ;
%Straight Lines [id:da2366819759458958] 
\draw [color={rgb, 255:red, 74; green, 74; blue, 74 }  ,draw opacity=1 ] [dash pattern={on 4.5pt off 4.5pt}]  (335.33,102.47) -- (394.2,102.47) ;
%Shape: Rectangle [id:dp033646750912228685] 
\draw  [color={rgb, 255:red, 208; green, 2; blue, 27 }  ,draw opacity=1 ][fill={rgb, 255:red, 208; green, 2; blue, 27 }  ,fill opacity=0.3 ] (316.65,154.38) -- (393.95,102.27) -- (393.52,162.09) -- (316.22,214.2) -- cycle ;
%Shape: Right Triangle [id:dp562025169984061] 
\draw  [color={rgb, 255:red, 208; green, 2; blue, 27 }  ,draw opacity=1 ][fill={rgb, 255:red, 208; green, 2; blue, 27 }  ,fill opacity=0.3 ][line width=0.75]  (316.43,154.95) -- (257.33,213.8) -- (316.43,213.8) -- cycle ;

% Text Node
\draw (223.33,206.57) node [anchor=north west][inner sep=0.75pt]    {${\textstyle x}$};
% Text Node
\draw (465.83,165.17) node [anchor=north west][inner sep=0.75pt]    {$t$};
% Text Node
\draw (316.83,17.67) node [anchor=north west][inner sep=0.75pt]    {$\tau $};
% Text Node
\draw (315.33,93.67) node [anchor=north west][inner sep=0.75pt]    {$T$};
% Text Node
\draw (389.33,167.67) node [anchor=north west][inner sep=0.75pt]    {$T$};

\end{tikzpicture}
\caption{Prism $P$  \label{fig1}}
\end{center}
\end{figure}

The function $\hbox {w}$ satisfies the following relation:
$$
\hbox {w}_t+\hbox {w}_{\tau}-\sum_{i=1}^n a_i(\varepsilon,x,u_{x_i}(t,x))\hbox {w}_{x_ix_i}=
$$
$$
\sum_{i=1}^n\big[a_i(\varepsilon,x,u_{x_i}(t,x))-a_i(\varepsilon,x,u_{x_i}(\tau,x)\big]u_{x_ix_i}(\tau,x)+
$$
$$
\sum_{i=1}^n\big[b_i(\varepsilon,x,u_{x_i}(t,x))-b_i(\varepsilon,x,u_{x_i}(\tau,x)\big]+
$$
$$
+F(x,u(t,x),\nabla u(t,x))-F(x,u(\tau,x),\nabla u(\tau,x)).
$$
Introduce the function 
$$
\omega \equiv \hbox {w}\, e^{-\tau}
$$
which satisfies in $P$ the following linear ultraparabolic equation
$$
L\,\omega \equiv \omega_t+\omega_{\tau}+\omega-\sum_{i=1}^na_i(\varepsilon,x,u_{x_i}(t,x))\omega_{x_ix_i}=
$$
$$
e^{-\tau}\sum_{i=1}^n\big[a_i(\varepsilon,x,u_{x_i}(t,x))-a_i(\varepsilon,x,u_{x_i}(\tau,x)\big]u_{x_ix_i}(\tau,x)+
$$
$$
e^{-\tau}\sum_{i=1}^n\big[b_i(\varepsilon,x,u_{x_i}(t,x))-b_i(\varepsilon,x,u_{x_i}(\tau,x)\big]+
$$
$$
e^{-\tau}[F(x,u(t,x),\nabla u(t,x))-F(x,u(\tau,x),\nabla u(\tau,x))]. \leqno (3.10)
$$
Let
$$
\Gamma_{\tau}=
\partial P\setminus \{(t,\tau,x):t=T,\,0<\tau<T,\,x\in \Omega\}. \leqno(3.11)
$$
Suppose that the function $\omega$ attains its positive maximum at the point $N\in {\overline P}\setminus \Gamma_{\tau}$. At this point it should be 
$$
L\,\omega > 0,
$$
since $-\omega_{x_ix_i}(N)\ge 0,$ $\omega_t(N)\ge 0,$ $\omega_{\tau}(N)\ge0,$ $\omega(N)>0$. On the other hand 
at this point $\nabla \omega=0$ i.e. 
$$
u(t,x)>u(\tau,x) \ \ \hbox{and}\ \ \nabla u(t,x)=\nabla u(\tau,x)
$$ 
and hence, from (3.10) and (3.2),
$$
L\,\omega\Big|_{N}\leq 0.
$$
From this contradiction we conclude that $\omega$ can not attain its positive maximum in 
${\overline P}\setminus \Gamma_{\tau}$. 

Consider $\omega$ on $\Gamma_{\tau}$: 

\par\noindent
for $x\in \partial\Omega$, $t\in[0,T]$, $\tau\in[0,T]$ we have $\omega=-K(t-\tau)e^{-\tau}\leq 0$;
\par\noindent
for $t=\tau$, $x\in {\overline \Omega}$, $t\in[0,T]$ we have $\omega=0$;
\par\noindent
for $\tau=0,$ $t\in[0,T],$ $x\in {\overline \Omega}$ we have $\omega=u(t,x)-u_0(x)-K t\leq 0$ due to Lemma 3.2.

Consequently $\omega\leq 0$ in ${\overline P}$ i.e.
$$
u(t,x)-u(\tau,x)\leq K(t-\tau). \leqno (3.12)
$$

\medskip

Now subtracting (3.7) from (3.8) for ${\tilde v}(t,\tau,x)\equiv u(\tau,x)-u(t,x)$ we obtain 
$$
{\tilde v}_t+{\tilde v}_{\tau}-\sum_{i=1}^n a_i(\varepsilon,x,u_{x_i}(\tau,x)){\tilde v}_{x_ix_i}=
$$
$$
\sum_{i=1}^n\big[a_i(\varepsilon,x,u_{x_i}(\tau,x))-a_i(\varepsilon,x,u_{x_i}(t,x)\big]u_{x_ix_i}(t,x)+
$$
$$
\sum_{i=1}^n\big[b_i(\varepsilon,x,u_{x_i}(\tau,x))-b_i(\varepsilon,x,u_{x_i}(t,x)\big]+
$$
$$
F(x,u(\tau,x),\nabla u(\tau,x))-F(x,u(t,x),\nabla u(t,x)).
$$
As in the previous case we conclude that the function
$$
{\tilde \omega} \equiv \tilde {\hbox {w}}\, e^{-\tau}\equiv \big({\tilde v}-K(t-\tau)\big)e^{-\tau}
$$
satisfies in $P$ the following ultraparabolic equation
$$
{\tilde \omega}_t+{\tilde \omega}_{\tau}+{\tilde \omega}-\sum_{i=1}^na_i(\varepsilon,x,u_{x_i}(\tau,x)){\tilde \omega}_{x_ix_i}=
$$
$$
e^{-\tau}\sum_{i=1}^n\big[a_i(\varepsilon,x,u_{x_i}(\tau,x))-a_i(\varepsilon,x,u_{x_i}(t,x)\big]u_{x_ix_i}(t,x)+
$$
$$
e^{-\tau}\sum_{i=1}^n\big[b_i(\varepsilon,x,u_{x_i}(\tau,x))-b_i(\varepsilon,x,u_{x_i}(t,x)\big]+
$$
$$
e^{-\tau}[F(x,u(\tau,x),\nabla u(\tau,x))-F(x,u(t,x),\nabla u(t,x))]. \leqno (3.13)
$$
Similarly to the previous case we obtain that ${\tilde \omega}$ can not attain its positive maximum in 
${\overline P}\setminus \Gamma_{\tau}$ and that ${\tilde \omega}\leq 0$ on $\Gamma_{\tau}$. The only difference is that 
for $\tau=0,$ $t\in[0,T],$ $x\in{\overline \Omega}$ we have ${\tilde \omega}=u_0(x)-u(t,x)-K t$ which is also 
non positive due to Lemma 3.2.

Consequently ${\tilde \omega}\leq 0$ in ${\overline P}$ i.e.
$$
u(\tau,x)-u(t,x)\leq K(t-\tau). \leqno (3.14)
$$
From (3.12) and (3.14) we conclude that in ${\overline P}$ 
$$
|u(t,x)-u(\tau,x)|\leq K(t-\tau).
$$
\smallskip
Taking into account the symmetry of the variables $t$ and $\tau$, we similarly consider the case $t<\tau$ and  obtain that 
in
$$
\{(t,\tau,x):t\in [0,T],\,\tau\in[0,T],x\,\in {\overline \Omega}\}
$$
the inequality
$$
|u(t,x)-u(\tau,x)|\leq K|t-\tau|
$$
holds. The last implies the required estimate. $\square$

\medskip

2. \underline{Isotropic case}.

Consider equation (2.4) coupled with conditions (1.2), (1.3). Denote by $K_1$ the following quantity:
$$
K_1=\max\Big|div \Big((|\nabla u_0|^2+\varepsilon)^{p(x)/2}\nabla u_0\Big)\Big| +\max|F(x,u,\nabla u_0)|
$$ 
where maximum is taken over the set $(\varepsilon,x,u)\in [0,\varepsilon_0]\times 
\overline{\Omega}\times [-M,M]$. One can easily see that $K_1<+\infty$ due to (1.10). 

\smallskip

\par\noindent
{\bf Lemma 3.4}. {\it For every $\varepsilon\in (0,\varepsilon_0]$ 
the following estimate
$$
|u_{\varepsilon}(t,x)-u_0(x)|\leq K_1\,t,
$$
takes place.
}
\par\noindent 
{\bf Proof.} Introduce the function
$$
h(t)=(K_1+\delta)t\ \ \hbox{in}\ \ [0,T],
$$
where $\delta>0$. Let us prove the following inequality
$$
u(t,x)-u_0(x)\leq h(t). \leqno (3.15)
$$
Consider the linear operator 
$$
L\equiv \frac{\partial}{\partial t}-\sum_{ij=1}^na_{ij}(\varepsilon,x,\nabla u_0)\frac{\partial^2}{\partial x_i\partial x_j}.
$$
Define the function $\phi^+\equiv u-[u_0(x)+h(t)]$, obviously
$$
L\phi^+=u_t-\sum_{ij=1}^n a_{ij}(\varepsilon, x, \nabla u_0)u_{x_ix_j}
+\sum_{ij=1}^n a_{ij}(\varepsilon, x, \nabla u_0)u_{0x_ix_j}-K_1-\delta.
$$
At the point $N\in {\overline Q_T}\setminus\Gamma_T$  (see (3.5) ) of possible maximum of the 
function $\phi^+$ we have 
$$
\nabla \phi^+=0\ \ \ \Leftrightarrow \ \ \nabla u=\nabla u_0 \ \ \Rightarrow\ \ 
a_{ij}(\varepsilon,x, \nabla u_0)=
a_{ij}(\varepsilon,x,\nabla u)
$$
and similarly to the anisotropic case we obtain 
$$
L\phi^+\Big|_{N}=div \Big((|\nabla u_0|^2+\varepsilon)^{p(x)/2}\nabla u_0\Big)+F(x,u,\nabla u_0)-K_1-\delta\Big|_{N}<0
$$
which is impossible. Obviously (see anisotropic case) $\phi^+\leq 0$ on $\Gamma_T$ and consequently $\phi^+\leq 0$ in $Q_T$. Thus (3.15) is proved.

\medskip

Let us show now that
$$
u(t,x)-u_0(x)\ge -h(t).\leqno (3.16)
$$
Again, similarly to the previous, at the point 
$N_1\in {\overline Q_T}\setminus\Gamma_T$ of possible minimum of the function $\phi^-\equiv u-[u_0(x)-h(t)]$ we have
$$
\nabla \phi^-=0\ \ \Leftrightarrow\ \ \nabla u=\nabla u_0\ \ \Rightarrow\ \ a_{ij}(\varepsilon,x,\nabla u_0)=
a_{ij}(\varepsilon,x,\nabla u)\ \ \Rightarrow
$$
$$
L\phi^-\Big|_{N_1}=div \Big((|\nabla u_0|^2+\varepsilon)^{p(x)/2}\nabla u_0\Big)+F(x,u,\nabla u_0)+K+\delta\Big|_{N_1}>0,
$$
which is impossible. Taking into account that $\phi^-\ge 0$ on $\Gamma_T$ we 
obtain the needed estimate (3.16) which with (3.15) imply 
$$
|u(t,x)-u_0(x)|\leq h(t).
$$
Passing to the limit when $\delta\to 0$ we finish the prove of Lemma 3.4. $\square$

\smallskip

\noindent
{\bf Lemma 3.5}. {\it For every $\varepsilon\in (0,\varepsilon_0]$ 
the inequality 
$$
|u_{\varepsilon t}|\leq K_1
$$
holds.}
\par\noindent 
{\bf Proof.} Consider equation (2.4) in two different points $(t,x)$ and $(\tau, x)$:
$$
u_t-\sum_{ij=1}^n a_{ij}(\varepsilon,x,\nabla u)u_{x_ix_j}=\sum_{i=1}^nb_i(\varepsilon,x,\nabla u)+F(x,u,\nabla u),\ \ u=u(t,x) \leqno(3.17)
$$
$$
u_{\tau}-\sum_{ij=1}^n a_{ij}(\varepsilon,x,\nabla u)u_{x_ix_j}=\sum_{i=1}^nb_i(\varepsilon,x,\nabla u)+F(x,u,\nabla u),\ \ u=u(\tau,x). \leqno(3.18)
$$
Subtracting (3.18) from (3.17) for $v(t,\tau,x)\equiv u(t,x)-u(\tau,x)$, since 
$$
v_t=u_t(t,x),\ \ v_{\tau}=-u_{\tau}(\tau,x),\ \ v_{x_ix_j}= u_{x_ix_j}(t,x)-u_{x_ix_j}(\tau,x),
$$
we obtain 
$$
v_t+v_{\tau}-\sum_{ij=1}^n a_{ij}(\varepsilon,x,\nabla u(t,x))v_{x_ix_i}=
$$
$$
\sum_{ij=1}^n\big[a_{ij}(\varepsilon,x,\nabla u(t,x)-a_{ij}(\varepsilon,x,\nabla u(\tau,x))\big]u_{x_ix_i}(\tau,x)+
$$
$$
\sum_{i=1}^n[b_i(\varepsilon, x,\nabla u(t,x))-b_i(\varepsilon, x,\nabla u(\tau,x))]+
F(x,u(t,x),\nabla u(t,x))-F(x,u(\tau,x),\nabla u(\tau,x)).
$$
Similarly to the anisotropic case, the function
$$
\omega \equiv \big(u(t,x)-u(\tau,x)-K_1(t-\tau)\big) e^{-\tau}
$$
satisfies in the domain $P$ (see 3.9) the following linear ultraparabolic equation
$$
L\,\omega \equiv \omega_t+\omega_{\tau}+\omega-\sum_{ij=1}^na_{ij}(\varepsilon,x,\nabla u(t,x))\omega_{x_ix_j}=
$$
$$
e^{-\tau}\sum_{ij=1}^n\big[a_{ij}(\varepsilon,x,\nabla u(t,x)-a_i(\varepsilon,x,\nabla u(\tau,x))\big]u_{x_ix_j}(\tau,x)
+
$$
$$
e^{-\tau}\sum_{i=1}^n[b_i(\varepsilon, x,\nabla u(t,x))-b_i(\varepsilon, x,\nabla u(\tau,x))]+
$$
$$
e^{-\tau}[F(x,u(t,x),\nabla u(t,x))-F(x,u(\tau,x),\nabla u(\tau,x))]. \leqno (3.19)
$$
Suppose that the function $\omega$ attains its positive maximum at the point $N\in {\overline P}\setminus \Gamma_{\tau}$ (see (3.11)). At this point it should be $L\,\omega >0,$
since matrix $a_{ij}$ is positively defined, matrix $\omega_{x_ix_j}(N)$ non positively defined  
and $\omega_t(N)\ge 0,$ $\omega_{\tau}(N)\ge 0,$ $\omega(N)>0$. On the other hand 
at this point $\nabla \omega=0$ i.e. 
$$
u(t,x)>u(\tau,x) \ \ \hbox{and}\ \ \nabla u(t,x)=\nabla u(\tau,x)
$$ 
and hence, from (3.19) and (3.2),
$$
L\,\omega\Big|_{N}\le 0.
$$
From this contradiction we conclude that $\omega$ can not attain its positive maximum in 
${\overline P}\setminus \Gamma_{\tau}$. 

Consider $\omega$ on $\Gamma_{\tau}$: 

\par\noindent
for $x\in \partial\Omega$, $t\in[0,T]$, $\tau\in[0,T]$ we have $\omega=-K(t-\tau)e^{-\tau}\leq 0$; for $t=\tau$, $x\in {\overline \Omega}$, $t\in[0,T]$ we have $\omega=0$; for $\tau=0,$ $t\in[0,T],$ 
$x\in {\overline \Omega}$ we have $\omega=u(t,x)-u_0(x)-K t\leq 0$ due to Lemma 3.4.

Consequently $\omega\leq 0$ in ${\overline P}$ i.e.
$$
u(t,x)-u(\tau,x)\leq K_1(t-\tau). \leqno (3.20)
$$

\medskip

Now subtracting (3.17) from (3.18) for ${\tilde v}(t,\tau,x)\equiv u(\tau,x)-u(t,x)$ we obtain 
$$
{\tilde v}_t+{\tilde v}_{\tau}-\sum_{ij=1}^n a_{ij}(\varepsilon,x,\nabla u(\tau,x)){\tilde v}_{x_ix_i}=
$$
$$
\sum_{ij=1}^n\big[a_{ij}(\varepsilon,x,\nabla u(\tau,x)-a_{ij}(\varepsilon,x,\nabla u(t,x))\big]u_{x_ix_j}(t,x)+
$$
$$
\sum_{i=1}^n\big[b_i(\varepsilon,x,\nabla u(\tau,x))-b_i(\varepsilon,x,\nabla u(t,x))\big]+
$$
$$
+F(x,u(\tau,x),\nabla u(\tau,x))-F(x,u(t,x),\nabla u(t,x)).
$$
The function 
$$
{\tilde \omega} \equiv \big(\tilde {v}-K_1(t-\tau)\big)\, e^{-\tau}
$$
in $P$ satisfies the following ultraparabolic equation
$$
{\tilde \omega}_t+{\tilde \omega}_{\tau}+{\tilde \omega}-\sum_{i=1}^na_{ij}(\varepsilon,x,u_{x_i}(\tau,x)){\tilde \omega}_{x_ix_i}=
$$
$$
e^{-\tau}\sum_{ij=1}^n\big[a_{ij}(\varepsilon,x,\nabla u(\tau,x)-a_i(\varepsilon,x,\nabla u(t,x))\big]u_{x_ix_j}(\tau,x)+
$$
$$
e^{-\tau}\sum_{i=1}^n\big[b_i(\varepsilon,x,\nabla u(\tau,x))-b_i(\varepsilon,x,\nabla u(t,x))\big]+
$$
$$
e^{-\tau}[F(x,u(\tau,x),\nabla u(\tau,x))-F(x,u(t,x),\nabla u(t,x)).
$$
Similarly to the previous case we obtain that ${\tilde \omega}$ can not attain its positive maximum in 
${\overline P}\setminus \Gamma_{\tau}$ and that ${\tilde \omega}\leq 0$ on $\Gamma_{\tau}$. The only difference is that for $\tau=0,$ $t\in[0,T],$ $x\in{\overline \Omega}$ we have ${\tilde \omega}=u_0(x)-u(t,x)-K t$ which is also non positive due to Lemma 3.4.

Consequently ${\tilde \omega}\leq 0$ in ${\overline P}$ and we have
$$
u(\tau,x)-u(t,x)\leq K(t-\tau). \leqno (3.21)
$$
From (3.20) and (3.21) we conclude that in ${\overline P}$ 
$$
|u(t,x)-u(\tau,x)|\leq K(t-\tau).
$$
\smallskip
Taking into account the symmetry of the variables $t$ and $\tau$ we similarly consider the case $t<\tau$ to obtain that in
$$
\{(t,\tau,x):t\in [0,T],\,\tau\in[0,T],x\,\in {\overline \Omega}\}
$$
the inequality
$$
|u(t,x)-u(\tau,x)|\leq K|t-\tau|
$$
holds. The last implies the required estimate. $\square$

\smallskip

{\bf Remark.} Concerning the linear and nonlinear ultraparabolic equations see [1], [19] and the references therein.

\smallskip

\bigskip

\centerline {{\bf \S 4. A priori estimates of $u_{\varepsilon x_i}$}}

\medskip

1. \underline {Anisotropic case}.

Here as in the previous section we take more general in compare with (1.4) right hand side of the equation. Namely we suppose that
$$
F=\sum_{i=1}^{n}f_i(x,u,u_{x_i})+f_0(x,u), \leqno (4.1)
$$
where for $i=1,...,n$ 
$$
f_i(x,u,u_{x_i})\leq \alpha|u_{x_i}|^{\tilde{p}_i+2}+\beta  
$$
with arbitrary positive constants $\alpha, \beta$ and constant $\tilde{p}_i$ such that 
$$
-1<\tilde{p}_i< \min_{\overline{\Omega}} p_i(x). 
$$

\noindent
{\bf Lemma 4.1}. {\it For the solution of problem $(2.1)$, $(1.2)$, $(1.3)$ there are constants 
$\tilde{C}_i$ independent of $\varepsilon$ such that  
$$
\max_{[0,T]}\int_{\Omega}|u_{\varepsilon x_i}(x,t)|^{p_i(x)+2}\,dx\le C_i, \ \ 
i=1,\dots,n. 
$$
}
{\bf Proof.} Multiplying equation (2.1) by $u_t$ and integrating by parts we obtain
$$
\int_{\Omega}u_t^2dx+
\sum_{i=1}^n\frac{1}{p_i+2}\int_{\Omega}\Big((u_{x_i}^{2}+\varepsilon)^{\frac{p_i+2}{2}}\Big)_tdx =
\int_{\Omega} Fu_t dx\leq K\int_{\Omega} |F|dx.
$$
Integrating with respect to $t$ we have 
$$
\sum_{i=1}^n\frac{1}{p_i+2}\int_{\Omega}(u_{x_i}^{2}+\varepsilon)^{\frac{p_i+2}{2}}dx
\leq K\int_{Q_T} |F|dx+\sum_{i=1}^n\frac{1}{p_i+2}\int_{\Omega}(u_{0x_i}^{2}+\varepsilon)^{\frac{p_i+2}{2}}dx.
$$
Taking into account that $g(\varepsilon)\equiv (z_i^{2}+\varepsilon)^{p_i+2/2}$ is an increasing function we obtain:
$$
\sum_{i=1}^n\frac{1}{p_i+2}\int_{\Omega}|u_{x_i}|^{p_i+2}dx 
\leq K\int_{Q_T} \sum_{i=0}^n|f_i|dx+\tilde{C}_1\leq 
$$
$$
\int_0^TK\Big(\alpha\sum_{i=1}^n\int_{\Omega}|u_{x_i}|^{\tilde{p_i}+2}dx+(n\beta+
\max_{\overline{\Omega}\times[-M,M]}|f_0|)|\Omega|\Big)dt+\tilde{C}_1,
$$
where
$$
\tilde{C}_1=\sum_{i=1}^n\frac{1}{p_i+2}\int_{\Omega}(u_{0x_i}^{2}+\varepsilon_0)^{\frac{p_i+2}{2}}dx.
$$
Now, taking into account condition (4.1) and using the H\"{o}lder inequality with non constant exponents 
(see [12], Theorem 2.1) and then the Young inequality  we have
$$
\int_{\Omega}|u_{x_i}|^{\tilde{p_i}+2}dx\leq c_p|\Omega|^{\frac{p_i-\tilde{p}_i}{p_i+2}}
\Big(\int_{\Omega}|u_{x_i}|^{p_i+2}dx\Big)^{\frac{\tilde{p}_i+2}{p_i+2}}\leq
c_p\frac{p_i-\tilde{p}_i}{p_i+2}|\Omega|+c_p\frac{\tilde{p}_i+2}{p_i+2}\int_{\Omega}|u_{x_i}|^{p_i+2}dx,
$$
constant $c_p$ in our case is equal to $|\Omega|+1$ (see [12])). Thus for 
$$
y(t)=\sum_{i=1}^n\frac{1}{p_i+2}\int_{\Omega}|u_{x_i}|^{p_i+2}dx 
$$
we have 
$$
y(t)\leq \tilde{C}_1+\int _0^T\big(\tilde{C}_2+\tilde{C}_3y(t)\big)dt. \leqno (4.2)
$$
Here 
$$
\tilde{C}_2=K\big(n\alpha (p_0+1)c_p|\Omega|+
(n\beta+\max_{\overline{\Omega}\times [-M,M]} |f_0(x,u)|)|\Omega|\big),\ \ 
\tilde{C}_3=c_pK\alpha(p_0+2),
$$
here $p_0=\max\{\max_{\Omega}p_1,..., \max_{\Omega}p_n\}$. 
From (4.2) by Gronwall's lemma we obtain the claimed result. $\square$

\medskip

2. \underline {Isotropic case}.

Here instead of (4.1)we suppose that 
$$
F=f(x,u,\nabla u)+f_0(x,u), \leqno (4.3)
$$
where
$$
f(x,u,\nabla u)\leq \alpha|\nabla u|^{\tilde{p}+2}+\beta 
$$
with arbitrary positive constants $\alpha, \beta$ and constant $\tilde{p}$ such that 
$$
-1<\tilde{p}<\min_{\overline{\Omega}} p(x). 
$$ 

\noindent
{\bf Lemma 4.2}. {\it For the solution of problem $(2.3), (1.2), (1.2)$  there exists a constant $C$ 
independent of $\varepsilon$ such that  
$$
\max_{[0,T]}\int_{\Omega}|\nabla u|^{p+2}\,dx\le C, \ \ 
i=1,\dots,n.
$$
}
{\bf Proof.} Multiplying equation (2.3) by $u_t$ and integrating first by parts over $x$, and then
integrating over $t$, similarly to the previous case we obtain 
$$
\frac{1}{p+2}\int_{\Omega}|\nabla u|^{p+2}dx
\leq K\int_{Q_T} |F|dx+\frac{1}{p+2}\int_{\Omega}\big(|\nabla u_0|^2+\varepsilon_0\big)^{\frac{p+2}{2}}dx\leq
$$
$$
\int_0^TK\Big(\alpha\int_{\Omega}|\nabla u|^{\tilde{p}+2}dx+(\beta+\max|f_0|)|\Omega|\Big)dt
+\frac{1}{p+2}\int_{\Omega}\big(|\nabla u_0|^2+\varepsilon_0\big)^{\frac{p+2}{2}}dx.
$$
Just as in the previous lemma, we have  
$$
\int_{\Omega}|\nabla u|^{\tilde{p}+2}dx\leq 
c_p|\Omega|^{\frac{p-\tilde{p}}{p+2}}\Big(\int_{\Omega}|\nabla u|^{p+2}dx\Big)^{\frac{\tilde{p}+2}{p+2}}
\leq
c_p\frac{p-\tilde{p}}{p+2}|\Omega|+c_p\frac{\tilde{p}+2}{p+2}\int_{\Omega}|\nabla u|^{p+2}dx. 
$$
Thus for 
$$
y(t)=\frac{1}{p+2}\int_{\Omega}|\nabla u|^{p+2}dx 
$$
we have 
$$
y(t)\leq \tilde{C}_1+\int _0^T\big(\tilde{C}_2+\tilde{C}_3y(t)\big)dt
$$
whence, by Gronwall's lemma, we obtain the required estimate. Here 
$$
\tilde{C}_1=\frac{1}{p+2}\int_{\Omega}\big(|\nabla u_0|^2+\varepsilon_0\big)^{\frac{p+2}{2}}dx,
$$
$$
\tilde{C}_2=K\big(\alpha c_p (p_0+1)|\Omega|+(\beta+\max|f_0|)|\Omega|\big),\ \ 
\tilde{C}_3=K\alpha c_p(p_0+2),\ \ p_0=\max p(x).
$$

\bigskip

\smallskip

\centerline {{\bf \S 5. Proof of Theorems 1, 2}}

\medskip

We obtain a weak solution to problem (1.1) - (1.3)
as a limit of the approximate solutions $u_\varepsilon$ constructed
in the previous section. Multiplying equation (2.1) by an arbitrary smooth function $\phi$ which vanishes on $(0,T)\times \partial\Omega$ 
and integrating by parts we have
$$
\int_{Q_T} u_{\varepsilon t}\,\phi\,dxdt+\int_{Q_T} \sum_{i=1}^n (u_{\varepsilon x_i}^2+\varepsilon)^{p_i/2}\, u_{\varepsilon x_i}\,\phi_{x_i} \,dxdt =
\int_{Q_T} \sum_{i=1}^n\big[f_i(x)u_{\varepsilon x_i}+f_0(x,u_{\varepsilon})\big]\phi dxdt. \leqno (5.1)
$$
As it follows from Lemmas 3.3, 4.1,  there exists a sequence $\varepsilon_k$ such that
$$
u_{\varepsilon_k x_i}\to u_{x_i} \ \ \hbox{weakly in}\ \ L^{\tilde{p}_i+2}(\Omega)\ \ \hbox{for}\ \ 
\tilde{p}_i=\min_{\overline{\Omega}}p_i(x),
$$
$$
u_{\varepsilon_k t}\to u_t \ \ \hbox{*-weakly in}\ \ L^\infty(Q_T)
$$
as $\varepsilon_k\to 0$. From the well known compactness Lemma (see [15, Chapter 1, \S 5 ], [8]) the estimates 
of Lemmas 3.3, 4.1 imply 
$$
u_{\varepsilon_kx_i}\to u_{x_i} \ \  \hbox{in} \ \ L^r(Q_T)\ \ \hbox{norm for some}\ \  r\in (1,+\infty).
$$
If, in particular, $\mu=\min_i \tilde{p}_i + 2\ge n$, then $r\in (1,+\infty)$ is arbitrary, otherwise 
$r\in(1, n\mu/(n-\mu)$. Thus, in order to pass to the limit in (5.1), we only have to prove that
$$
\int_{Q_T} \sum_{i=1}^n (u_{\varepsilon_k x_i}^2
+\varepsilon_k)^{p_i/2}\, u_{\varepsilon_k x_i}\,\phi_{x_i}\, dxdt\to
\int_{Q_T} \sum_{i=1}^n |u_{x_i}|^{p_i}
\, u_{x_i}\,\phi_{x_i}\, dxdt \ \ \hbox{as}\ \ \varepsilon_k\to 0. 
$$
This can be done similarly to [18, p.3018] or [24, p.175].

\smallskip

The uniqueness of the weak solution can be proved by standard considerations taking into account the monotonicity of the elliptic part of the operator (see [18, p. 3020).  

Theorem~1 is proved. The proof of Theorem~2 is similar.

\medskip

{\bf Remark about more general right hand side.} Note that the estimates obtained is Section 3 
and 4 for problem (2.3), (1.2), (1.3) are valid for 
$$
F=f(x,u,\nabla u)+f_0(x,u),  \leqno (5.2)
$$ 
satisfying only the monotonicity condition (3.2) and growth restrictions (4.3). The obtained estimates are not enough to pass to the limit in $f(x,u_{\varepsilon},\nabla u_{\varepsilon})$ and obtain weak solution in the sense of Definition 2. Therefore, we take $F$ as in (1.4). However, if we assume that $p+2>n$ then $W^1_{q+2}$, with $q=\min_{\overline{\Omega}} p(x)$, is compactly embedded to $\bf{C}^{\alpha}$ (for some $\alpha\in (0,1)$) and thus $u_{\varepsilon}\to u$ uniformly. This allows us to study the above problem in the framework of viscosity solutions with $F$ as in (5.2) with constraints (3.2), (4.3) (see, for example, [23]). 

Similarly for the problem (2.1), (1.2), (1.3) the estimates in Sections 3, 4 are valid for 
$$
F=\sum_{i=1}^{n}f_i(x,u,u_{x_i})+f_0(x,u), 
$$ 
satisfying only the monotonicity condition (3.2) and growth restrictions (4.1). Thus, if $p_i+2>n$ for all $i$, then, similarly to the previous case, we obtain that $u_{\varepsilon}\to u$ uniformly and can  study the formulated problem in the framework of viscosity solutions .

\medskip

\bigskip

\centerline{{\bf REFERENCES}}

\smallskip

\par\noindent 
[1] F. Anceschi, S.Polidoro, \emph{A survey on the classical theory for Kolmogorov equation}, Matematiche 75, No. 1, 221-258 (2020). 

\par\noindent 
[2] S. Antontsev, S. Shmarev, \emph{Anisotropic parabolic equations with variable nonlinearity}, Publ. Mat. 53 (2), 355--399 (2009).

\par\noindent 
[3] S. Antontsev, S. Shmarev, \emph{Vanishing solutions of anisotropic parabolic equations with variable nonlinearity}, J. Math. Anal. Appl. 361, No. 2, 371-391 (2010). 

\par\noindent 
[4] S. Antontsev, I. Kuznetsov, S. Shmarev, \emph{Global higher regularity of solutions to singular p(x, t)-parabolic equations}, J. Math. Anal. Appl. 466 (1) (2018) 238--263.

\par\noindent 
[5] S. Antontsev, S. Shmarev, \emph{Higher regularity of solutions of singular parabolic equations with variable nonlinearity}, Appl. Anal. 98 (1 \& 2), 310--331 (2019).

\par\noindent 
[6] S. Antontsev, S. Shmarev, \emph{Global estimates for solutions of singular parabolic and elliptic equations with variable nonlinearity}, Nonlinear Anal. 195, 111--724, (2020).

\par\noindent 
[7] S. Antontsev, V. Zhikov, \emph{Higher integrability for parabolic equations of p(x, t)-Laplacian type, Adv. Differ. Equ.} 10 (9), 1053--1080  (2005).

\par\noindent
[8] J.P. Aubin, \emph{Un theoreme de compacite}, C.R. acad. sci. Paris, 256, 5042 - 5044, (1963). 

\par\noindent 
[9] M.M. Boureanu, A. Velez-Santiago, \emph{Fine regularity for elliptic and parabolic anisotropic Robin problems with variable exponents}, J. Differ. Equations 266, No. 12, 8164-8232 (2019). 

\par\noindent 
[10] V. B\"{o}gelein, F. Duzaar, \emph{H\"{o}lder estimates for parabolic $p(x,t)$-Laplacian systems}, Math. Ann. 354 (3), 907--938  (2012).

\par\noindent 
[11] A.H. Erhardt, \emph{Compact embedding for $p(x,t)$-Sobolev spaces and existence theory to parabolic equations with $p(x,t)$-growth}, Rev. Mat. Complut. 30 (1), 35--61 (2017).

\par\noindent
[12] O. Kovacik, J. Rakosnik, \emph{On spaces $L_{p(x)}$ and $W_{k,p(x)}$} 
Czech. Math. J. 41(116), No. 4, 592-618 (1991).
 
\par\noindent 
[13] O.A.Ladyzhenskaya, V.A.Solonnikov, N.N.Ural'tseva, \emph{Linear and Quasilinear Equations of Parabolic Type.} American Mathematical Society, Providence R.I.,1968, 648pp. 

\par\noindent 
[14] S. Lian, W. Gao, H. Youan, C. Cao, \emph{Existence of solutions to an initial Dirichlet problem of evolutional $p(x)$-Laplace equations}, Ann. Inst. Henri Poincare, Anal. Non Lineaire 29, No. 3, 377-399 (2012). 

\par\noindent 
[15] Lions, J.-L., \emph{Quelques methodes de resolution des problemes aux limites non lineaires}, (1969), Dunod, Gauthier-Villars Paris, p. 554 

\par\noindent 
[16] S. Shmarev, \emph{On the continuity of solutions of the nonhomogeneous evolution $p(x,t)$-Laplace equation}, Nonlinear Anal. 167, 67--84  (2018). 

\par\noindent 
[17] J. Simsen, M. Simsen, M. Teixeira Primo, \emph{On $p_s(x)$-Laplacian parabolic problems with non-globally Lipschitz forcing term}, Z. Anal. Anwend. 33(4), 447-462 (2014).

\par\noindent 
[18] V.N. Starovoitov, Al.Tersenov, \emph{Singular and degenerate anisotropic parabolic equations with a nonlinear source}, Nonlinear Anal. 72 , no. 6, 3009--3027 (2010).

\par\noindent 
[19] Al. Tersenov, \emph{On the global solvability of the Cauchy problem for a quasilinear ultraparabolic equation}, Asymptotic Anal. 82, No. 3-4, 295-314 (2013).

\par\noindent 
[20] Al. Tersenov, \emph{The one dimensional parabolic $p(x)$
-Laplace equation}, NoDEA, Nonlinear Differ. Equ. Appl. 23, No. 3, Paper No. 27, 11 p. (2016). 

\par\noindent 
[21] Al. Tersenov, Ar. Tersenov, \emph{On the Bernstein-Nagumo's condition in the theory of nonlinear parabolic equations}, J. Reine Angew. Math. 572, 197-217 (2004). 

\par\noindent
[22] Al. Tersenov, Ar. Tersenov, \emph{Existence of Lipschitz continuous solutions to the Cauchy-Dirichlet problem for anisotropic parabolic equations}, J. Funct. Anal. 272, No. 10, 3965-3986 (2017). 

\par\noindent 
[23] Al. Tersenov, Ar. Tersenov,
\emph{Existence results for anisotropic quasilinear parabolic equations with time-dependent exponents and gradient term}, J. Math. Anal. Appl. 480, No. 1, Article ID 123386, 18 p. (2019). 

\par\noindent 
[24] Al. Tersenov, Ar. Tersenov, \emph{On quasilinear anisotropic parabolic equations with time-dependent exponents}, Sib. Math. J. 61, No. 1, 159-177 (2020). 

\par\noindent 
[25] M. Xu, Y.Z. Chen, \emph{H\"{o}lder continuity of weak solutions for parabolic equations with nonstandard growth conditions}, Acta Math. Sin. Engl. Ser. 22 (3), 793--806 (2006). 
 
\par\noindent  
[26] F. Yao, \emph{H\"{o}lder regularity for the general parabolic $p(x,t)$-Laplacian equations}, NoDEA Nonlinear Differ. Equ. Appl. 22 (1), 105--119 (2015).

\par\noindent 
[27] H. Zhan, Z. Feng, \emph{Existence and stability of the doubly nonlinear anisotropic parabolic equation}, 
J. Math. Anal. Appl. 497, No. 1, Article ID 124850, 23 p. (2021). 

\par\noindent 
[28] H. Zhan, Z. Feng, \emph{Solutions of evolutionary $p(x)$-Laplacian equation based on the weighted variable exponent space}, Z. Angew. Math. Phys. 68, No. 6, Paper No. 134, 17 p. (2017).  

\par\noindent 
[29] C. Zhang, S. Zhou, X. Xue, \emph{Global gradient estimates for the parabolic $p(x,t)$-Laplacian equation}, Nonlinear Anal. 105, 86--101  (2014).

\end{document}